\documentclass[12pt]{article}
\usepackage{theorem, amsfonts,amsmath,graphicx}
\usepackage{setspace}        

\nonstopmode
\newcommand\qed{\hfill$\sqcap\kern-8.0pt\hbox{$\sqcup$}$}

\title{Getting real with Real Options}

\author{M. R. Grasselli \\ 
McMaster University}

\begin{document}
\maketitle

\begin{abstract}
We apply a utility--based method to obtain the value of a finite--time investment opportunity when the underlying real asset is not perfectly correlated to a traded financial asset. Using a discrete--time algorithm to calculate the indifference price for this type of real option, we present numerical examples for the corresponding investment thresholds, in particular highlighting their dependence with respect to correlation and risk aversion.  
\end{abstract}


\section{Introduction}
Most of the standard literature in real options is based on one or both of the following unrealistic assumptions: (1) that the time horizon for the problem at hand is infinite and (2) that the real asset under consideration is perfectly correlated to a traded financial asset. The {\em infinite--maturity} hypothesis helps to reduce the dimensionality of the problem by removing its dependence on time, therefore concentrating on stationary solutions only, while the {\em spanning asset} hypothesis allows the introduction of useful replication arguments developed for derivative pricing in complete markets. Together they led to the development of a coherent and intuitive approach for investment under uncertainty, well--represented for instance in Dixit and Pindyck (1994), where the decisions to start, abandon, reactive and mothball a given project were reduced to the solution of systems of linear equations. 

Since then, several authors have dropped the artifice of an infinite maturity time, therefore resourcing to numerical methods for dealing with the non-stationary valuation problem. These include finite--difference methods for the associated partial differential equation and lattice methods for discrete--time option pricing. However, even recent books such as Smit and Trigeorgis (2004) carry the assumption that ``real-options valuation is still applicable provided we can find a reliable estimate for the market value of the asset" (page 102), 
which is tantamount to saying that ``markets are sufficiently complete".  In reality, most investment problems where the real options approach is deemed relevant occur in markets which are far from being complete. For example, almost by definition an R\&D investment decision concerns a product which is {\em not} currently commercialized and therefore commands uncertainty that is at best imperfectly correlated with available financial assets. 

Exceptions to this adherence to a ``near completion" assumption, but still in the context of an infinite time horizon, are Hugonnier and Morellec (2004)\nocite{HugoMore04} and Henderson (2005)\nocite{Hend05}. In the first paper, a risk--averse manager facing an investment decision tries to maximize his expected utility considering the effect that shareholders' external control will have on his personal wealth. By assuming that the underlying project is subject to both market risk, which the manager can hedge using a traded financial asset, and idiosyncratic risk, which cannot be hedged in the available financial market, the authors reduce this decision to an investment problem in an incomplete market. By contrast, under a similar model for a project with both systematic and idiosyncratic risks, Henderson (2005) uses an exponential utility framework in order to actually calculate the value for the opportunity to invest as a derivative in an incomplete market, therefore remaining closer in spirit to the real options paradigm. 

In this paper, we study a finite--horizon version of Henderson's model. We first review the mechanism for pricing derivatives in an incomplete markets using an exponential utility in the context of investment decisions in a simple one--period binomial model in section \ref{one-period-section}. In section \ref{multi-period-section} we extend the valuation procedure to a multi--period model used an an approximation for continuous--time markets, followed by numerical experiments exploring the properties of the option to invest in 
section \ref{numerical}, including comparisons with the corresponding infinite horizon and complete markets limits. Section \ref{conclusion-section} then presents conclusions drawn from the model, especially in contrast with alternative ways of dealing with market incompleteness in the context of real options.    

\section{Investment decisions in one period}                                                                                                                                                                                                                                                                                                                                                         
\label{one-period-section}

Consider an investor who needs to decide whether to pay a sunk cost $I$ in exchange of a project with current value $V_0$. Assume  that such investment can be made either at time $0$ or postponed until time $T$, when the project value is expected to rise or fall according to specified probabilities, so that the opportunity to invest is formally equivalent to discrete-time early-exercise call option with strike price $I$ having the project value as the underlying asset. When the project value is perfectly correlated to the price of a traded financial asset, such option can be priced using standard arbitrage and replication arguments. In the absence of such spanning asset, the option to invest becomes analogous to a derivative in an incomplete market. Instead of wishfully pretending that risk--neutral and replication arguments can still be used in this case, we argue that the investor's risk preference should be explicitly used for valuing the option to invest. For this, we follow Hobson and Henderson (2002) \nocite{HendHobs02} and consider a utility indifference framework based on an exponential utility of the form $U(x)=-e^{-\gamma x}$.
 
Let us assume the existence of a riskless cash account with constant annualized  interest rate $r$, which we use as a fixed {\em numeraire}, and denote the discounted project value by $V$ and the discounted price of a correlated traded financial asset by $S$. We then specify their one--period dynamics by 
\begin{equation}
(S_T,V_T)=\left\{\begin{array}{l}
(uS_0,hV_0) \quad \mbox{ with probability } p_1, \\
(uS_0,\ell V_0) \quad \mbox{  with probability } p_2, \\
(dS_0,h V_0) \quad \mbox{ with probability } p_3, \\
(dS_0,\ell V_0) \quad \mbox{  with probability } p_4, \end{array}\right.
\label{one-period}
\end{equation}
where $0<d<1<u$ and $0<\ell<1 <h$, for positive initial values $S_0,V_0$ and historical probabilities $p_1,p_2,p_3,p_4$.

Without the opportunity to invest in the project $V$, a rational agent with initial wealth $x$ will keep an amount $\xi$ in the cash account  and holding $H$ units of the traded asset $S$ in such a way as to maximize the expected utility of the terminal wealth
\begin{equation}
X^x_T=\xi+HS_T=x+H(S_T-S_0).
\end{equation}
That is, the investor will try to solve the optimization problem 
\begin{equation}
\label{Merton}
u^0(x)=\max_H E[U(X^x_T)]
\end{equation}

Suppose next that the investor is offered the opportunity to invest in the project at the end of the period, thereby receiving a discounted payoff $C_T=(V_T-e^{-rT}I)^+$. In order to acquire this option, the investor must spend an amount of money $\pi$ at time $t=0$. For instance, $\pi$ might be the price of land that will allow a subsequent real estate development, or the price of a license to explore a natural resource. The key assumption now is that, after obtaining the option, the investor will try to use the financial market in such a way as to maximize the utility of the total terminal position. In other words, an investor with initial wealth $x$ who acquires the option for the price $\pi$ will try to solve the modified optimization problem
\begin{equation}
u^C(x-\pi)=\sup_{h} E[U(X^{x-\pi}_T+C_T)]
\label{optimalhedge}
\end{equation}

Following Hodges and Neuberger (1989)\nocite{HodNeu89}, we define the {\em indifference price} for the option to invest in the final period as the amount $\pi^C$ that solves the equation 
\begin{equation}
u^0(x)=u^C(x-\pi).
\label{ind_id}
\end{equation}

Denoting the two possible pay-offs at the terminal time by $C_h$ and $C_\ell$, it is then a straightforward calculation to show that, for an exponential utility, such indifference price is given by 
\begin{equation}
\label{european}
\pi^C=g(C_h,C_\ell)
\end{equation}
where, for fixed parameters $(u,d,p_1,p_2,p_3,p_4)$ the function $g:\mathbb{R}\times\mathbb{R}\rightarrow \mathbb{R}$ is given by 
\begin{equation}
g(x_1,x_2)=\frac{q}{\gamma}\log \left(\frac{p_1+p_2}{p_1e^{-\gamma x_1}+p_2e^{-\gamma x_2}}\right)+
\frac{1-q}{\gamma}\log\left(\frac{p_3+p_4}{p_3e^{-\gamma x_1}+p_4e^{-\gamma x_2}}\right),
\label{gfunction}
\end{equation}
with
\[q=\frac{1-d}{u-d}.\]

We will henceforth refer to $\pi^C$ as the {\em continuation value} for holding the option of investing at a later time. If we now introduce the possibility of investment at time $t=0$, it is clear that immediate exercise of this option will occur whenever its {\em exercise value} $(V_0-I)^+$ is larger than its continuation value $\pi^C$. That is, from the point of view of this agent, the value at time zero for the opportunity to invest in the project either at $t=0$ or $t=T$ is given by 
\[C_0=\max\{(V_0-I)^+,g((hV_0-e^{-rT}I)^+,(\ell V_0-e^{-rT}I)^+)\}.\]

\section{The multi-period model}
\label{multi-period-section}

Consider now a continuous-time two--factor market of the form 
\begin{eqnarray}
dS_t &=& (\mu_1-r)S_tdt +\sigma_1 S_t dW \\
dV_t &=& (\mu_2-r)V_tdt +\sigma_2 V_t (\rho dW + \sqrt{1-\rho^2}) dZ,
\end{eqnarray}
where $\mu_1,\mu_2\sigma_1,\sigma_2,r$ are constants and $(W,Z)$ are standard independent Brownian motions under the historical probability measure $P$. As before, $S_t$ corresponds to the discounted price of a traded financial asset which is correlated to the discounted value of the underlying project $V_t$ with a correlation parameter $\rho$.  

An approximation for this market can be obtained by dividing the time interval $[0,T]$ into $N$ subintervals with equal time steps 
$\Delta t = T/N$ and taking the one--period dynamics for the discrete--time processes $(S_n,V_n)$ to be given by \eqref{one-period}.  We then need to choose the dynamic parameters $u,d,h,\ell$ and the one-period probabilities $p_i$ so that, in the limit of small $\Delta t$, such dynamics matches the distributional properties of the continuos time processes $S_t$ and $V_t$. For instance, one can verify that, up to terms of order $\Delta t$, the choices 
\begin{eqnarray}
u &=& e^{\sigma_1\sqrt{\Delta t}}, \qquad h=e^{\sigma_2\sqrt{\Delta t}}, \\
d &=& e^{-\sigma_1\sqrt{\Delta t}}, \qquad \ell=e^{-\sigma_2\sqrt{\Delta t}},\\
p_1+p_2&=&\frac{e^{(\mu_1-r)\Delta t}-d}{u-d}, \\
p_1+p_3&=&\frac{e^{(\mu_2-r)\Delta t}-\ell}{h-\ell} \\
\rho \sigma_1 \sigma_2\Delta t &=& (u-d)(h-\ell)[p_1p_4-p_2p_3],
\end{eqnarray}
provide the correct variances, means and correlation for $S_t$ and $V_t$. Supplemented by the condition 
\begin{equation}
p_1+p_2+p_3+p_4=1,
\end{equation}
these equations uniquely determine the historical probabilities $p_i$. 

Having fixed these parameters, let us choose a sufficiently large integer $M$ and denote
\begin{equation}
V^{(i)}=h^{M+1-i}V_0, \qquad i=1,\ldots, 2M+1,
\label{array}
\end{equation}
ranging from $(h^MV_0)$  to $(\ell^M V_0)$, respectively the highest and lowest achievable discounted project values 
starting from the middle point $V_0$ with the multiplicative parameter $h=\ell^{-1}>1$. In practice, $M$ should be chosen so that such highest and lowest values are comfortably beyond the range of project values that can be reached during the time interval $[0,T]$ with reasonable probabilities (say by corresponding to returns which are away from their mean by more than four standard deviations). 
Then each realization for the discrete-time process $V_n$ following the dynamics \eqref{one-period} can then be thought of as a path over a $(2M+1)\times N$ rectangular grid having the values \eqref{array} as its repeated columns.  

The discounted value of the option to invest on the project can then be determined as a function $C_{in}$ on this grid, with the index $i=1,\ldots, 2N+1$ referring to the underlying project value $V^{(i)}$, and the index $n=0,\ldots,N$ referring to time 
$t_n=n\Delta t$.  We start by specifying the following boundary conditions:
\begin{eqnarray}
C_{iN}&=&(V^{(i)}-e^{-rT}I)^+, \qquad i=1,\ldots,2N+1, \label{terminal} \\
C_{1n}&=& V^{(1)}-e^{-rn\Delta_t}, \qquad \quad n=0,\ldots,N, \label{top} \\
C_{2N+1,n}&=& 0, \qquad \qquad \qquad \qquad n=0,\ldots,N. \label{bottom}
\end{eqnarray}
Condition \eqref{terminal} corresponds to the fact that at maturity the option to invest should be exercised whenever the project value exceeds the investment cost, whereas conditions \eqref{top} and \eqref{bottom} mean that such option should always be exercised when the project value is at its highest but is worthless when the project is at its lowest. The values in the interior of the grid are then obtained by backward induction as follows:
 \begin{equation}
C_{in}=\max\left\{(V^{(i)}-e^{-rn\Delta t}I)^+, g(C_{i+1,n+1},C_{i-1,n+1})\right\}, \quad \begin{array}{l}
n=N-1,\ldots, 0 \\ i=2,\ldots, 2N. \end{array}
\end{equation}
That is, at each node on the grid, the investor chooses between exercising the investment option, obtaining its immediate exercise value 
$V^{(i)}-e^{-rn\Delta t}I$ or holding the option one step into the future, retaining its continuation value $g(C_{i+1,n+1},C_{i-1,n+1})$. 
 
Accordingly, at each time $t_n$, the exercise threshold $V^*_n$ is defined as the project value for which the exercise value for the option becomes higher than its continuation value. For project values below $V^*_n$, the investor will prefer to hold  the option, while for project values higher than such threshold, preference for immediate exercise will prevail.  

\section{Numerical Experiments}
\label{numerical}

We now investigate how the exercise threshold, and consequently the value of the option to invest, varies with the different model parameters. In what follows, unless explicitly indicated, we fix the investment cost $I=1$, the risk--free interest $r=0.04$, the time--to--maturity $T=10$, the dynamic parameters for the traded asset $\mu_1=0.115$, $\sigma_1=0.25$, $S_0=1$  and the volatility for the project 
$\sigma_2=0.2$. Given these parameters, the CAPM equilibrium expected rate of return on the project for a given correlation 
$\rho$ is
\begin{equation}
\bar\mu_2=r+\rho\left(\frac{\mu_1-r}{\sigma_1}\right)\sigma_2.
\label{CAPM}
\end{equation}
The difference $\delta=\bar\mu_2-\mu_2$, known as the below--equilibrium rate--of--return shortfall, should then be interpreted as the incomplete market analogue of a dividend rate paid by the project, which we fix at $\delta=0.04$. 

\begin{figure}[h]
\begin{center}
\includegraphics[width=\textwidth]{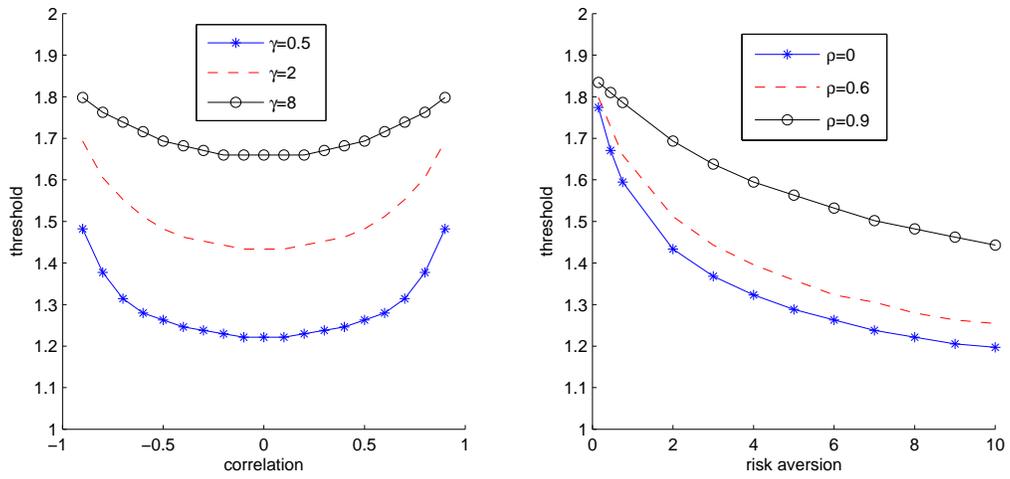}
\end{center}
\caption{Exercise threshold as a function of correlation and risk aversion.}
\label{corr_gamma}
\end{figure}

Because incompleteness is the main theme of this paper, we starting with the dependence with correlation. In the complete market case, corresponding to $\rho\rightarrow 1$ in our model, the investment threshold can be obtained in closed-form in the limit of an infinite time horizon. For the parameters above, such formula gives 
$V_{\scriptscriptstyle DP}^*=2$ (see page 153 in Dixit and Pindyck (1994)). By contrast, the investment threshold obtained from a simple NPV criterion (that is, invest whenever the net present value for the project is positive) in this case is equal to 
$V^*_{\scriptscriptstyle NPV}=1$. This constitutes the most widespread result from real option theory: irreversibility and time flexibility 
lead to investors waiting until much larger thresholds before committing to an investment decision. In the left panel of figure \ref{corr_gamma} we see what reservations need to be made in the presence of market incompleteness. Predictably, the exercise 
thresholds decrease as we move away from $\rho=\pm 1$, meaning that the presence of unhedged risks decreases the value of the option to invest. Interestingly however, even at its minimum, corresponding to $\rho=0$, the investment threshold is still higher than suggested by NPV. That is, even when the risk in the project is entire idiosyncratic and therefore cannot be hedged with financial assets, time flexibility still confers an option value to the opportunity to invest which is higher than its net present value, {\em irrespectively of any replication argument}.     

\vspace{0.2in}
\centerline{[Insert figure \ref{corr_gamma} here]}
\vspace{0.2in}

Moving to the dependence with risk--aversion, for an incomplete market and still in the infinite time horizon setup, McDonald and Siegel (1986)\nocite{McDSie86} obtain the exercise threshold by assuming that investors require compensation for market risks whilst being risk--neutral towards idiosyncratic risk. As observed in Henderson (2005), this corresponds to the limit $\gamma\rightarrow 0$ in an exponential utility framework. For the parameters above, with the adjustment for market risks done according to CAPM, this threshold coincides with the value $V_{\scriptscriptstyle}^*=2$ of Dixit and Pindyck (1994). We can observe these limiting behavior in the right panel of figure \ref{corr_gamma}, together with the fact that the investment threshold for our model expectedly decreases as a function of the parameter $\gamma$. That is, risk--aversion can significantly erode the option value obtained from time flexibility. In the limit $\gamma\rightarrow\infty$, one can explicitly show that expression \eqref{gfunction} tends to the subhedge price of the derivative, which is zero for a call option, so that value for the investment opportunity reduces to its net present value, with the corresponding investment 
threshold $V^*_{\scriptscriptstyle NPV}=1$. As we observe in the graph, this erosion of value with risk aversion is faster the lower the correlation between the project and the traded asset.

\begin{figure}[h]
\begin{center}
\includegraphics[width=\textwidth]{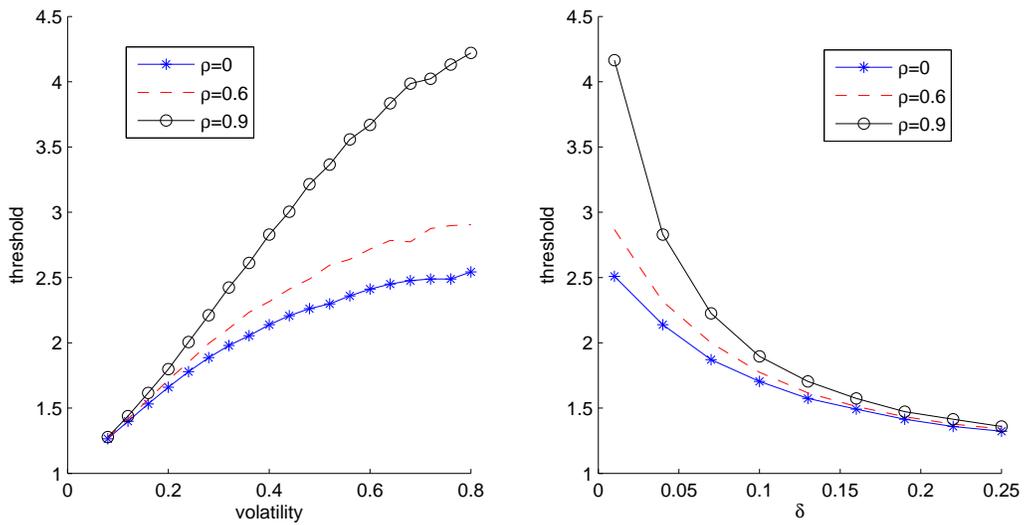}
\end{center}
\caption{Exercise threshold as a function of volatility and dividend rate.}
\label{vol_delta}
\end{figure}

As with classical real options, we can verify in our model that heightened uncertainty in project values increases the value of time flexibility,  
which is reflected in the increasing threshold with respect to the volatility $\sigma_2$ observed in the left panel of figure \ref{vol_delta}. In the same figure, we observe how the threshold decreases as a function of $\delta$, implying that the incentives for immediate investment in the project are higher the more its rate of return falls below its equilibrium value, which is the incomplete market analogue of more investment when dividends are higher.

\vspace{0.2in}
\centerline{[Insert figure \ref{vol_delta} here]}
\vspace{0.2in}

\begin{figure}[h]
\begin{center}
\includegraphics[width=\textwidth]{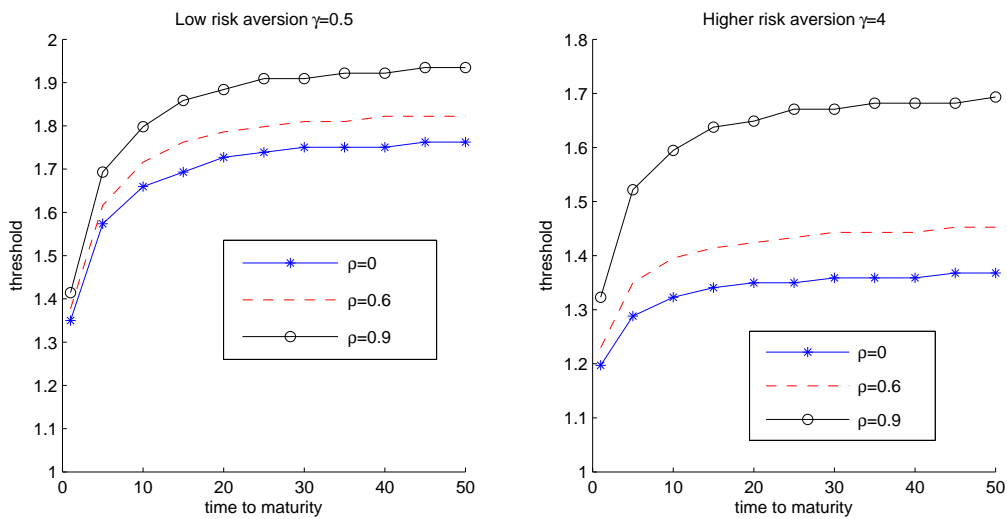}
\end{center}
\caption{Exercise threshold as a function of time to maturity.}
\label{maturity}
\end{figure}

All of these features were already present in the elegant closed-form solution obtained in Henderson (2005) for the infinite time horizon limit and under the somewhat artificial restriction that the investment cost should grow at a rate $\alpha=r$. Our contribution consists on the one hand in removing this restriction (our investment cost can grow at any rate, in particular at a rate $\alpha=0$ as above) and on the other hand in numerically calculating the exercise threshold as a function of time--to--maturity. In figure \ref{maturity} we see that the exercise threshold can take a long time to converge to its asymptotic value, particularly in the desirable cases of low risk--aversion and high--correlation, indicating that for typical maturities of only a couple of years the stationary solution provides a poor approximation for its finite--horizon counterpart.    

\vspace{0.2in}
\centerline{[Insert figure \ref{maturity} here]}
\vspace{0.2in}

\begin{figure}[h]
\begin{center}
\includegraphics[width=\textwidth]{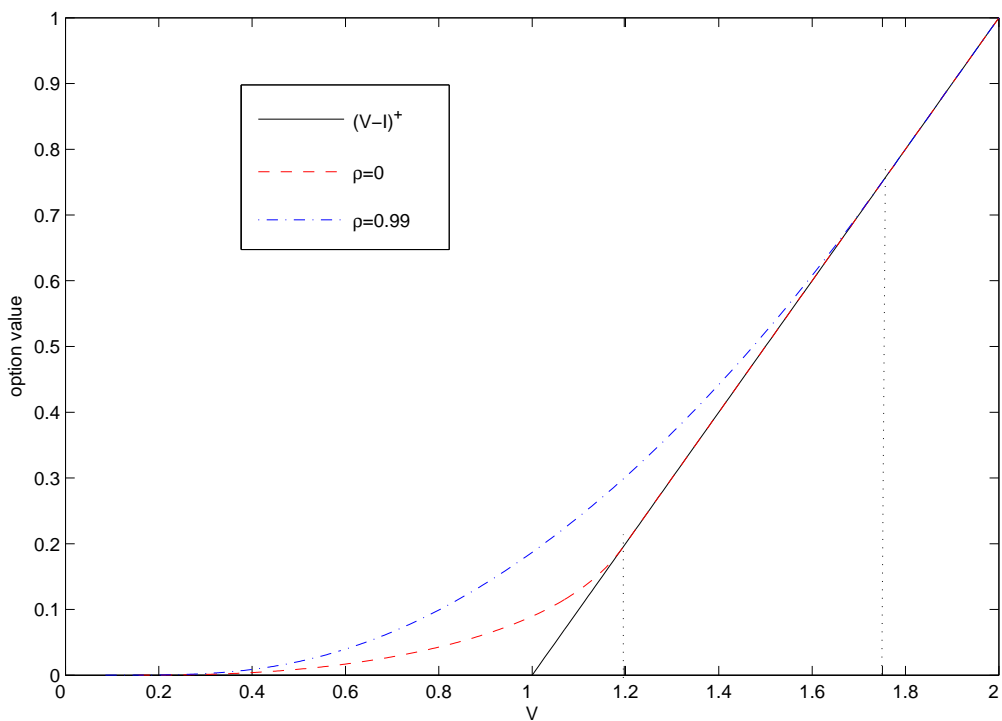}
\end{center}
\caption{Option value as a function of underlying project value. The threshold for $\rho=0$ is 1.1972 and the one for $\rho=0.99$ is 1.7507.}
\label{value}
\end{figure}

We conclude this section with a graph of option values as a function of the current level of the underlying project. These should be compared with either figure 5.3 on page 154 of Dixit and Pindyck (1994) or figure 2 of Henderson (2005), which use the same base parameters as ours, except for time--to--maturity, which we take to be $T=10$. We can observe once more that higher correlations lead to higher option values. Moreover, we confirm our previous observation that even for $\rho=0$ the opportunity to invest is more valuable then its net present value, represented in the graph by the solid line depicting the function $(V-I)^+$. Finally, notice how the smooth pasting and matching conditions, which were not {\em a priori} assumed in our model, are satisfied by the option values, in the sense that the curves in figure \ref{value} smoothly match the function $(V-I)^+$ at the corresponding exercise thresholds, marked in the graph by the two vertical dotted lines.  

In all of the numerical experiments above, we used a fixed time step $\Delta t=1/900$, so that the relative precision for project values on the grid is of the order $\sigma_2\sqrt{\Delta t}\sim 0.0067$. For each point marked in the pictures above, the thresholds and option values were obtained on a typical $1000\times 9000$ grid in approximately 15 seconds using a desktop PC at 3MHz. 

\section{Discussion}
\label{conclusion-section}

We have proposed a multi--period binomial model for assessing the value of an option to invest on a project in a finite time horizon in the absence of a perfectly spanning financial asset. The exercise thresholds obtained from our model exhibit the expected dependence with respect to correlation, uncertainty, risk aversion, dividend rates and time to maturity. In particular, for perfectly correlated or perfectly anti--correlated assets the threshold approaches that of a complete market, whereas when the aversion to idiosyncratic risks tends to zero the threshold approaches the one obtained under similar assumptions by McDonald and Siegel (1986). Moreover, we verify that even in the zero correlation case, whereby none of the risk in the project can be hedge in a financial market, the paradigm of real options can still be applied to value an investment decision, since the opportunity to invest still carries an option value above its net present value. In other words, it is time flexibility itself, more than the possibility of replication, that is the source of the extra value of an investment opportunity. This value, however, erodes sharply at higher levels of risk aversion, and even more so when the project is uncorrelated to financial markets. 

Apart from the outright use of risk--neutral valuation even when markets are incomplete - under the wishful assumption that they are complete enough for all practical purposes - the most widespread alternative method for dealing with incompleteness in a real options context is through the use of dynamic programming with an exogenous discount rate. This is the approach indicated for instance in the second half of chapter 5 in Dixit and Pindyck (1994), in which an investor equates the expected capital appreciation from a project to the expected rate of return  on the investment opportunity, using a corporate rate of return, which is different from the risk--free interest rate and meant to express corporate risk preferences. Despite its popularity, such approach has the serious theoretical drawback that the  fully nonlinear risk preferences of a corporation can hardly be expressed through a single discount factor. In fact, the majority of financial economics literature use an expected utility function {\em together} with an exogenous discount factor in order to model risk preferences. 

At a more practical level, this dynamic programming approach with a corporate discount rate obscures the most important aspects of real options, namely the intuition that can be gain when managerial decisions are treated as options. For example, under the option paradigm, investment on a multi-stage project is analogous to a portfolio of options, each having its own value and interacting in a complex manner towards the value of the whole project. Precisely because such analogies are completely lost in the dynamic programming approach, 
authors such as Dixit and Pindyck dropped it in the remaining of their book in favor of a contingent claim analysis, which then formally relies back on the complete market framework with a spanning asset hypothesis.    

In comparison, our proposed method handles incompleteness by explicitly introducing risk preferences in an economically sound utility--based framework for the realistic case of a finite time horizon, while retaining the computational complexity of a standard binomial valuation. The method can then be extended in order to provide incomplete market versions for all the standard managerial decisions treated as real options. For example, all the thresholds for investment, abandonment, suspension and reactivation of a project in an incomplete market  are obtained in detail in the forthcoming paper by Grasselli and Nitzan (2006). As an application in a different 
direction, a modification of the method proposed here was used in Grasselli (2005) for the valuation of employee stock options. In this way, our valuation method considerably enlarges the domain of applicability of both real options and indifference pricing.

\end{document}